\documentclass[11pt]{article}

\usepackage{amssymb,latexsym}
\usepackage{amsmath}
\usepackage{cite}
\usepackage{eucal}

\title{A Remark on the Deformation of GNS Representations of $^*$-Algebras} 

\author{\textbf{Stefan Waldmann\thanks{Stefan.Waldmann@physik.uni-freiburg.de}}  
  \\[0.5cm]
  Fakult{\"a}t f{\"u}r Physik\\ 
  Albert-Ludwigs-Universit{\"a}t Freiburg\\
  Hermann Herder Stra{\ss}e 3\\
  D 79104 Freiburg\\
  Germany
  }

\date{December 2000\\[0.5cm] FR-THEP 2000/18}

\newcommand{\im} {{\mathrm i}}
\newcommand{\eu} {{\mathrm e}}

\newcommand{\cc}[1]      {\overline{{#1}}}

\newcommand{\id}         {{\mathsf{id}}}

\newcommand{\ring}[1]    {{\mathsf{{#1}}}}
\newcommand{\SP} [1]     {{\left\langle{{#1}}\right\rangle}}

\newcommand{\starWeyl} {\mathbin{*_{\scriptscriptstyle\mathrm{Weyl}}}}
\newcommand{\starWick} {\mathbin{*_{\scriptscriptstyle\mathrm{Wick}}}}

\newtheorem{lemma} {Lemma} [section]
\newtheorem{proposition} [lemma] {Proposition}
\newtheorem{theorem} [lemma] {Theorem}

\newtheorem{example}[lemma] {Example}

\newenvironment{proof}{\small{\sc Proof:}}{{\hspace*{\fill} $\square$\\}}

\numberwithin{equation}{section}

\begin{document}

\maketitle

\begin{abstract}
Motivated by deformation quantization we investigate the algebraic GNS
construction of $^*$-representations of deformed $^*$-algebras over
ordered rings and compute their classical limit. The question if a
GNS representation can be deformed leads to the deformation of
positive linear functionals. Various physical examples from
deformation quantization like the Bargmann-Fock and the
Schr{\"o}dinger representation as well as KMS functionals are
discussed. 
\end{abstract}

%
% some introduction
%

\section{Introduction}

Deformation quantization as introduced in \cite{BFFLS78} has shown to
be not only a powerful quantization scheme of finite-dimensional
systems, see
\cite{DL83b,Fed94a,OMY91,NT95a,NT95b,Kont99,GR99,BCG97,Kon97b,BBG98}
for existence and classification results, but 
it also yields various applications to string theory and
non-commutative field theories, see
e.g.~\cite{Schomerus99,JS2000,JSW2000,BuWa2000c},  
index theory \cite{NT95a,NT95b,Fed96}, quantum field theory
\cite{Dit90,CatFeld2000,DF2000}, and quantization of systems with
constraints and phase space reduction \cite{Fed98a,BHW2000}.

While the deformed algebras yield models for the observables, a notion
of states by using positive linear functionals was introduced in
\cite{BW98a} and used to obtain an algebraic version of the GNS
construction of $^*$-representations. In various examples this concept
of $^*$-representations on pre-Hilbert spaces over the formal power
series was shown to be physically reasonable and mathematically
fruitful \cite{BNW99a,BNW98a,BNPW98,Wal2000a}.

In this article we continue the more systematic investigation of the
$^*$-representation theory of (deformed) $^*$-algebras as started in
\cite{BuWa99a,BuWa2000b,BuWa2000a,BuWa2000c}, where the notion of a
deformed $^*$-representation was introduced. Here the positivity of
the inner product of the pre-Hilbert space plays a crucial role. Since
the GNS representations are certainly the most important ones we shall
focus on this sub-class of $^*$-representations. To this end we
study deformations of positive functionals and relate the GNS
representations of the deformed and undeformed positive linear
functionals: the classical limit of the GNS representation of the
deformed functional gives back the GNS representation of the undeformed
one. This statement is in so far non-trivial as it involves a two-step
quotient procedure where in the intermediate step the quotients are
not the same. Here the notion of the classical limit of a pre-Hilbert
space in the sense of \cite[Sect.~8]{BuWa99a} comes in crucially.
It also gives an easy criterion whether a GNS representation can be
deformed: one only has to deform the corresponding positive linear
functional which always can be done in the case of deformation
quantization of symplectic manifolds \cite[Prop.~5.1]{BuWa2000a}. We
explain and illustrate this in several examples like the Bargmann-Fock
representation and Schr{\"o}dinger representation and point out the
importance of the notion of a positive deformation of a $^*$-algebra.

The paper is organized as follows: In Section~2 we recall some basic
facts on $^*$-algebras and the GNS construction. Section~3 contains
the definition and first properties of the deformation theory of
$^*$-algebras and their $^*$-representations. In Section~4 we prove
the main result and in Section~5 several examples are discussed.

%
% *-algebras and GNS construction
%

\section{$^*$-Algebras and the GNS construction}

Let us briefly recall the notion of $^*$-algebras over ordered
rings and the GNS construction of $^*$-representations out of a
positive linear functional, see \cite{BW98a,BW97b,BuWa99a} for more
details.

By $\ring R$ we denote an ordered ring where our main examples needed
for physical applications are $\mathbb R$ and 
$\mathbb R[[\lambda]]$. Then $\ring C = \ring R(\im)$ with 
$\im^2 = -1$ will be the generalization of the complex numbers. If
$\ring R$ is ordered then $\ring R[[\lambda]]$ is ordered, too, whence
this notion is well adapted to formal deformation theory. A
\emph{pre-Hilbert space} $\mathcal H$ over $\ring C$ is a 
$\ring C$-module with a positive definite Hermitian inner product,
i.e. a map 
$\SP{\cdot,\cdot} : \mathcal H \times \mathcal H \to \ring C$ which is
$\ring C$-linear in the second argument and satisfies
$\SP{\phi,\psi} = \cc{\SP{\psi, \phi}}$ and $\SP{\chi,\chi} > 0$ for
all $\phi, \psi, \chi \in \mathcal H$, $\chi \ne 0$. 
Note that this notion is meaningful since $\ring R$ is ordered. By 
$\mathcal B(\mathcal H)$ we denote the algebra of endomorphism of
$\mathcal H$ having an adjoint with respect to $\SP{\cdot, \cdot}$. A
\emph{$^*$-algebra} $\mathcal A$ over $\ring C$ is an associative
algebra over $\ring C$ with an involutive anti-linear
anti-automorphism called the $^*$-involution, which shall be denoted
by $A \mapsto A^*$. In particular, $\mathcal B(\mathcal H)$ is a
$^*$-algebra and $A^*$ is simply given by the adjoint of $A$. A
\emph{$^*$-representation} of $\mathcal A$ on $\mathcal H$ is a
$^*$-homomorphism $\pi: \mathcal A \to \mathcal B(\mathcal H)$.

A $\ring C$-linear functional $\omega: \mathcal A \to \ring C$ is
called \emph{positive} if for all $A \in \mathcal A$
\begin{equation}
    \label{POsFuncDef}
    \omega(A^*A) \ge 0.
\end{equation}
In this case $\omega$ satisfies a Cauchy-Schwarz inequality, see
e.g.~\cite[Lem.~5]{BW98a}, and the space
\begin{equation}
    \label{GelfandIdeal}
    \mathcal J_{\omega} 
    := \{ A \in \mathcal A \; | \; \omega(A^*A) = 0 \}
\end{equation}
turns out to be a left ideal in $\mathcal A$, the so-called
\emph{Gel'fand ideal}. Thus one has a $\mathcal A$-left module
structure on  
$\mathcal H_\omega = \mathcal A \big/ \mathcal J_\omega$ which is
traditionally denoted by
\begin{equation}
    \label{GNSRep}
    \pi_\omega(A)\psi_B = \psi_{AB},
\end{equation}
where $\psi_B \in \mathcal H$ denotes the equivalence class of 
$B \in \mathcal A$. Moreover, $\mathcal H_\omega$ is endowed with a
pre-Hilbert space structure by setting
\begin{equation}
    \label{GNSProd}
    \SP{\psi_A, \psi_B} = \omega(A^*B).
\end{equation}
Then $\pi_\omega$ is easily shown to be a $^*$-representation of
$\mathcal A$ on $\mathcal H_\omega$, the so-called GNS representation
corresponding to $\omega$. If $\mathcal A$ is unital then $\pi_\omega$
is cyclic with cyclic vector $\psi_1$ and we have 
$\omega(A) = \SP{\psi_1, \pi_\omega(A)\psi_1}$, a property which
characterizes $\pi_\omega$ up to unitary equivalence. Hence it will be
convenient to consider unital $^*$-algebras in the following. See
\cite{BuWa99a,BuWa2000b} for results in the non-unital case and
\cite{BW98a,BW97b,BNW98a,BNW99a,BNPW98,Wal2000a,BuWa2000c} for various
applications to deformation quantization.

%
% deformations of everything... 
%

\section{Deformations of $^*$-Algebras and $^*$-Representations}

Now consider a formal associative deformation 
$\boldsymbol{\mathcal A} = (\mathcal A[[\lambda]], \star)$ of a
$^*$-algebra $\mathcal A$ over $\ring C$ in the sense of Gerstenhaber
\cite{Ger64,GS88}, i.e. $\star$ is an associative product deforming
the product of $\mathcal A$ in higher orders of $\lambda$. Then this
deformation is called \emph{Hermitian} if the $^*$-involution of
$\mathcal A$ (extended $\ring C[[\lambda]]$-anti-linearly to 
$\mathcal A[[\lambda]]$) is still a $^*$-involution for $\star$, see
\cite[Sect.~3]{BuWa2000a}. In principle one can also deform the
$^*$-involution but in the physical applications this is not necessary
and even not wanted: in deformation quantization $\mathcal A$ will play
the role of the classical observable algebra and Hermitian elements in
$\mathcal A$ correspond to measurable physical quantities, the
observables. Thus this characterization should be preserved under
quantization. Here the formal parameter $\lambda$ corresponds to
Planck's constant $\hbar$ whenever the formal series converge.

Suppose now that $\boldsymbol{\mathcal H}$ is a pre-Hilbert space over
$\ring C[[\lambda]]$. Then we consider the space 
$\boldsymbol{\mathcal H}_0 := 
\{\boldsymbol{\phi} \in \boldsymbol{\mathcal H} \; | \; 
\SP{\boldsymbol{\phi},\boldsymbol{\phi}}|_{\lambda=0} = 0\}$.
It turns out that $\boldsymbol{\mathcal H}_0$ is a
$\ring C[[\lambda]]$-submodule and 
$\mathcal H := \mathfrak C\boldsymbol{\mathcal H} := 
\boldsymbol{\mathcal H} \big/ \boldsymbol{\mathcal H}_0$,
now viewed as a $\ring C$-module, is a pre-Hilbert space over 
$\ring C$. Indeed, if 
$\mathfrak C: \boldsymbol{\mathcal H} \to \mathcal H$
denotes the projection then the inner product is just given by 
$\SP{\mathfrak C \boldsymbol{\phi}, \mathfrak C \boldsymbol{\psi}}
= \SP{\boldsymbol{\phi}, \boldsymbol{\psi}}|_{\lambda=0}$, where
$\boldsymbol{\phi}, \boldsymbol{\psi} \in \boldsymbol{\mathcal H}$.
Finally, this procedure is functorial whence we shall refer to
$\mathfrak C$ as the classical limit functor \cite[Sect~8]{BuWa99a}. 
\begin{example}
    Let $\mathcal H$ be a pre-Hilbert space over $\ring C$ and
    consider the $\ring C[[\lambda]]$-module 
    $\boldsymbol{\mathcal H}$. By $\lambda$-linear extension of the
    inner product, $\boldsymbol{\mathcal H}$ becomes a pre-Hilbert
    space over $\ring C[[\lambda]]$. It follows immediately that
    $\boldsymbol{\mathcal H}_0 = \lambda \boldsymbol{\mathcal H} =
    \lambda \mathcal H[[\lambda]]$ and thus 
    $\mathfrak C \boldsymbol{\mathcal H} \cong \mathcal H$
    again. Moreover, we note that in this case
    \begin{equation}
        \label{BHseries}
        \mathcal B (\boldsymbol{\mathcal H})
        =
        \mathcal B(\mathcal H) [[\lambda]].
    \end{equation}
\end{example}

However, the classical limit functor is defined on all pre-Hilbert
spaces over $\ring C[[\lambda]]$ and, in fact, we shall need the more
general situation in the sequel.

If $\boldsymbol{\pi}$ is a $^*$-representation of
$\boldsymbol{\mathcal A}$ on some pre-Hilbert space
$\boldsymbol{\mathcal H}$ over $\ring C[[\lambda]]$ then we can also
define the classical limit 
$\pi = \mathfrak C\boldsymbol{\pi}$ of the representation by setting
$\pi(A) \mathfrak C\boldsymbol{\psi} 
= \mathfrak C (\boldsymbol{\pi}(A)\boldsymbol{\psi})$.
One obtains a $^*$-representation $\pi$ of $\mathcal A$ on
$\mathcal H = \mathfrak C \boldsymbol{\mathcal H}$.
Again, $\mathfrak C$ is functorial, see \cite[Sect~8]{BuWa2000a}.

This raises the following questions: For a given deformation
$\boldsymbol{\mathcal A}$ of $\mathcal A$, which $^*$-representations
of $\mathcal A$ occur as classical limits of $^*$-representations of
$\boldsymbol{\mathcal A}$, i.e. which $^*$-representations can be
quantized? Moreover, one wants to understand in how many ways (up to
unitary equivalence) this is possible and which properties are
preserved under the classical limit and quantization.

A more restricted problem is the following: given a
$^*$-representation $(\mathcal H, \pi)$ of $\mathcal A$ can we find a
$^*$-representation $\boldsymbol{\pi}$ of $\boldsymbol{\mathcal A}$ on
$\boldsymbol{\mathcal H} = \mathcal H[[\lambda]]$ with 
$\mathfrak C\boldsymbol{\pi} = \pi$? It turns out that this leads too
quickly into obstructions and is thus too restrictive:

Consider the polynomials $\mathcal A = \mathbb C[z,\cc z]$ in two
variables with the $^*$-involution $z^* = \cc z$. Moreover, let
$\mathcal H = \mathbb C$ with the standard inner product and 
$\pi(1) = \id$, $\pi(z) = 0 = \pi(\cc z)$. This determines a
$^*$-representation of $\mathcal A$ on $\mathcal H$. Now consider the
Wick product
\begin{equation}
    \label{WickProduct}
    f \starWick g 
    = \sum_{r=0}^\infty \frac{(2\lambda)^r}{r!}
    \frac{\partial^r f}{\partial z^r}
    \frac{\partial^r g}{\partial \cc z^r}
\end{equation}
for $f, g \in \mathcal A[[\lambda]]$ which gives a Hermitian
deformation, see e.g.~\cite{BW97a}. Then a simple and straightforward
computation shows the following lemma:
\begin{lemma}
    \label{NoWickLem}
    There is no $^*$-representation $\boldsymbol{\pi}$ of 
    $(\mathcal A[[\lambda]], \starWick)$ on $\mathbb C[[\lambda]]$
    deforming $\pi$.
\end{lemma}
\begin{proof}
    Suppose there would be such a $\boldsymbol{\pi}$. Then necessarily
    $\boldsymbol{\pi} (1) = \id$ and $\boldsymbol{\pi}(z)$ as well as 
    $\boldsymbol{\pi}(\cc z)$ would have to start in order
    $\lambda$. But this is incompatible with 
    $z \starWick \cc z - \cc z \starWick z = 2\lambda 1$.
\end{proof}

%
% deformation of GNS representations
%

\section{Deformation of GNS representations}

Since the problem of deforming a $^*$-representation 
$(\mathcal H, \pi)$ of $\mathcal A$ in full generality seems to be
rather difficult we shall now restrict to a more specific class of
representations, namely the GNS representations.

Hence let $\omega: \mathcal A \to \ring C$ be a positive linear
functional and let 
$\boldsymbol{\mathcal A} = (\mathcal A[[\lambda]], \star)$ be a
Hermitian deformation of $\mathcal A$. Moreover, let
$\boldsymbol{\omega}: 
\boldsymbol{\mathcal A} \to \ring C[[\lambda]]$
be a $\ring C[[\lambda]]$-linear functional which therefore can be
written as 
$\boldsymbol{\omega} 
= \sum_{r=0}^\infty \lambda^r \boldsymbol{\omega}_r$
with $\boldsymbol{\omega}_r: \mathcal A \to \ring C$. It is called a
\emph{deformation of $\omega$} if $\boldsymbol{\omega}$ is a positive
functional for $\boldsymbol{\mathcal A}$ and 
$\boldsymbol{\omega}_0 = \omega$.
In this situation we can compute the classical limit of the GNS
representation induced by $\boldsymbol{\omega}$:
\begin{theorem}
    \label{MainTheorem}
    Let $\mathcal A$ be a $^*$-algebra over $\ring C$ and
    $\boldsymbol{\mathcal A} = (\mathcal A[[\lambda]], \star)$ a
    Hermitian deformation of $\mathcal A$. If 
    $\boldsymbol{\omega}: 
    \boldsymbol{\mathcal A} \to \ring C[[\lambda]]$
    is a deformation of a positive linear functional
    $\omega = \boldsymbol{\omega}_0: \mathcal A \to \ring C$ then the
    classical limit of the GNS representation 
    $(\boldsymbol{\mathcal H_\omega}, \boldsymbol{\pi_\omega})$ is
    canonically unitarily equivalent to the GNS representation
    $(\mathcal H_\omega, \pi_\omega)$ by the intertwiner
    \begin{equation}
        \label{TheIntertwiner}
        U: \mathfrak C \boldsymbol{\mathcal H_\omega}
        \ni \mathfrak C \psi_{\boldsymbol{A}} \mapsto 
        \psi_{A} \in \mathcal H_{\omega},
    \end{equation}
    where 
    $\boldsymbol{A} = \sum_{r=0}^\infty \lambda^r \boldsymbol{A}_r 
    \in \boldsymbol{\mathcal A}$ with $A = \boldsymbol{A}_0$.
\end{theorem}
\begin{proof}
    It is straightforward to check that $U$ is well-defined, unitary,
    and an intertwiner between $\mathfrak C \boldsymbol{\pi_\omega}$
    and $\pi_\omega$.
\end{proof}

This result is in so far non-trivial as the relation between the two
Gel'fand ideals $\boldsymbol{\mathcal J_\omega}$ and 
$\mathcal J_\omega$ may be quite complicated. Typically, there is
\emph{no} isomorphism like 
$\boldsymbol{\mathcal J_\omega} \cong \mathcal J_\omega[[\lambda]]$
inside $\mathcal A[[\lambda]]$: Using the terminology of
\cite[Def.~30]{BHW2000} the Gel'fand ideal 
$\boldsymbol{\mathcal J_\omega}$ is \emph{not} a deformation of the
submodule $\mathcal J_\omega$ in general.

This can be understood heuristically as follows. From
$\boldsymbol{\omega} (A^* \star A) = 0$ one obtains usually
\emph{more} conditions than only the zeroth order condition
$\omega(A^*A) = 0$. Hence $\boldsymbol{\mathcal J_\omega}$ is
typically `smaller' than
$\mathcal J_\omega[[\lambda]]$. Thus the pre-Hilbert space
$\boldsymbol{\mathcal H_\omega}$ is `bigger' than 
$\mathcal H_\omega[[\lambda]]$. However, the `additional' vectors in
$\boldsymbol{\mathcal H_\omega}$ have inner products which vanish in
the classical limit $\lambda=0$. Thus 
$\mathfrak C\boldsymbol{\mathcal H_\omega} \cong \mathcal H_\omega$
becomes possible. Note however, that in general there is no canonical
way to compare $\boldsymbol{\mathcal H_\omega}$ with 
$\mathcal H_\omega[[\lambda]]$.

According to the above theorem the problem of deforming a
representation reduces to the deformation of a positive functional, at
least in the case of a GNS representation. Hence we shall investigate
which classically positive functionals can actually be deformed.

%
% Examples
%

\section{Examples}

A first example of a deformed GNS representation is given by
\emph{faithful} positive linear functionals. A positive linear
functional $\omega: \mathcal A \to \ring C$ is called faithful if
$\omega(A^*A) > 0$ for all $A \ne 0$. It follows 
$\mathcal J_\omega = \{0\}$ and hence 
$\mathcal H_\omega = \mathcal A$ as $\ring C$-modules. The GNS
representation $\pi_\omega$ is just the usual left-action. Then the
following statement is obvious:
\begin{lemma}
    \label{FaithLem}
    Let $\omega: \mathcal A \to \ring C$ be a faithful positive linear
    functional and let 
    $\boldsymbol{\omega} = \sum_{r=0}^\infty \lambda^r 
    \boldsymbol{\omega}_r: \boldsymbol{\mathcal A} 
    \to \ring C[[\lambda]]$
    be a real $\ring C[[\lambda]]$-linear functional 
    (i.e. $\boldsymbol{\omega} (A^*) = \cc{\boldsymbol{\omega} (A)}$)
    with $\boldsymbol{\omega}_0 = \omega$. Then 
    $\boldsymbol{\omega}$ is a faithful positive $\ring
    C[[\lambda]]$-linear functional of $\boldsymbol{\mathcal A}$ and
    thus a deformation of $\omega$.
\end{lemma}

Thus a faithful positive linear functional can always be deformed
whence the above theorem can be applied.

The physical importance of faithful positive functionals comes from
the fact that in deformation quantization these functionals correspond
to thermodynamical states. In particular the \emph{KMS functionals} on
symplectic manifolds are faithful, see
\cite{BFLS84,BRW98a,BRW99,Wal2000a} for further details.

On the other hand there are strong obstructions on the $^*$-algebra
$\mathcal A$ to have a faithful positive linear functional at all: in
particular, it follows immediately that the \emph{minimal ideal}
$\mathcal J_{\scriptscriptstyle\mathrm{min}} (\mathcal A)$, which is
the intersection of all Gel'fand ideals, has to be trivial. For a
detailed discussion and further consequences of this matter we refer
to \cite{BuWa2000b}.

Another important case is obtained if the deformation 
$\boldsymbol{\mathcal A}$ is a \emph{positive deformation} of
$\mathcal A$, i.e. if any positive linear functional 
$\omega: \mathcal A \to \ring C$ can be deformed
\cite[Def.~4.1]{BuWa2000a}. In this case we have the following result:
\begin{proposition}
    \label{PosDefProp}
    Let $\boldsymbol{\mathcal A}$ be a positive deformation of a
    unital $^*$-algebra $\mathcal A$ over $\ring C$. If 
    $(\mathcal H, \pi)$ is an orthogonal sum of cyclic (and hence GNS)
    $^*$-representations then $(\mathcal H, \pi)$ can be deformed into
    a $^*$-representation for $\boldsymbol{\mathcal A}$.
\end{proposition}
\begin{proof}
    This follows immediately from the fact that the classical limit
    functor preserves orthogonal sums.
\end{proof}

Unfortunately, in general there may be $^*$-representations of
$\mathcal A$ which are not orthogonal sums of cyclic ones. Even for
$C^*$-algebras one needs to complete a direct sum in order to obtain
all (non-degenerate) representations. However, in
physical applications of deformation quantization, where $\mathcal A$
plays the role of the observable algebra, this can typically
be assumed in order to obtain interesting representations. Here a
cyclic vector has the interpretation of the `vacuum vector' or `ground
state': It should be possible to obtain all vectors by applying the
observables to the ground state. Furthermore, it was shown in
\cite[Prop.~5.1]{BuWa2000a} that star products on symplectic manifolds
are automatically positive deformations. Thus the above situation
applies to these star products of deformation quantization.

Let us now discuss more concrete examples from deformation
quantization. First we consider the Wick star product $\starWick$ on
$\mathbb C^n$ which generalizes (\ref{WickProduct}) to 
$C^\infty (\mathbb C^n)[[\lambda]]$, see e.g.~\cite{BW97a}. It is
known that $\starWick$ is even a strongly positive deformation,
i.e. all classically positive linear functionals are positive for
$\starWick$ as well \cite[Lem.~4.4]{BuWa2000a}. In particular, the
$\delta$-functional at $0 \in \mathbb C^n$ is positive and yields
the Bargmann-Fock representation \cite{BW98a}. We shall discuss
this example now in detail to illustrate the above introduced
notions. For the classical and quantum Gel'fand ideals $\mathcal J$
and $\boldsymbol{\mathcal J}$ one obtains
\begin{equation}
    \label{WickGelfandIdeal}
    \mathcal J 
    = \{f \in C^\infty (\mathbb C^n) \; | \; f(0) = 0 \} 
    \quad \textrm{and}\quad
    \boldsymbol{\mathcal J} 
    = 
    \Big\{
    f \in C^\infty (\mathbb C^n)[[\lambda]] 
    \; \Big| \;  
    \forall I \ni \mathbb N^n: \;
    \frac{\partial^{|I|} f}{\partial \cc z^I} (0) = 0
    \Big\},
\end{equation}
respectively. Hence the representation spaces $\mathcal H$ and
$\boldsymbol{\mathcal H}$ are given by
\begin{equation}
    \label{FockSpace}
    \mathcal H = \mathbb C
    \quad
    \textrm{and}
    \quad
    \boldsymbol{\mathcal H} = 
    \left(\mathbb C[[\cc y^1, \ldots, \cc y ^n]]\right)[[\lambda]]
\end{equation}
with the usual $\mathbb C$-valued inner product for $\mathcal H$ and
\begin{equation}
    \label{FockInnerProd}
    \SP{\psi, \phi} =
    \sum_{r=0}^\infty \frac{(2\lambda)^r}{r!} \sum_{|I| = r}
    \cc{\frac{\partial^r \psi}{\partial \cc y^I} (0)} \;
    \frac{\partial^r \phi}{\partial \cc y^I} (0),
\end{equation}
for $\psi, \phi \in \boldsymbol{\mathcal H}$. The GNS representation
is classically given by $\pi(f) 1 = f(0) 1$, analogously to
Lemma~\ref{NoWickLem}. Quantum mechanically it is
the usual Bargmann-Fock representation, in particular
$\boldsymbol{\pi} (z^k) = 2\lambda \frac{\partial}{\partial \cc y^k}$
and $\boldsymbol{\pi} (\cc z^k) = \cc y^k$, see
\cite[Lem.~8]{BW98a}. Now we see that $\boldsymbol{\mathcal J}$ is
strictly smaller than $\mathcal J[[\lambda]]$ and
$\boldsymbol{\mathcal H}_0$ contains all vectors in
$\boldsymbol{\mathcal H}$ except the $\mathbb C$-multiples of
$1$. Hence it is \emph{not} just $\lambda \boldsymbol{\mathcal H}$.
Thus we see that with the more general notion of a deformation of a
$^*$-representation we can avoid the obstructions from
Lemma~\ref{NoWickLem}. Let us also remark that the same phenomena
occurs for the Wick star product on an arbitrary K{\"a}hler manifold
$M$ and the $\delta$-functional at some $p \in M$, cf.~\cite{BW98a}.
We also remark that in this example the classical limit $\pi$ is `more
degenerate' than the quantized representation $\boldsymbol{\pi}$. On
e.g. the polynomials on $\mathbb C^n$ the representation
$\boldsymbol{\pi}$ is faithful while $\pi$ is not. See
\cite[Prop.~8.5]{BuWa99a} for a more detailed discussion of which
properties are preserved under the classical limit.

In a last example we shall discuss the formal Schr{\"o}dinger
representation as discussed in \cite{BW98a,BNW98a,BNW99a,BNPW98}. We
consider the cotangent bundle 
$\pi: T^*\mathbb R^n \to \mathbb R^n$ with the usual Weyl-Moyal star
product 
\begin{equation}
    \label{WeylMoyal}
    f \starWeyl g = \mu \circ 
    \eu^{\frac{\im\lambda}{2} \sum_k \left(
          \frac{\partial}{\partial q^k} \otimes
          \frac{\partial}{\partial p_k} 
          -
          \frac{\partial}{\partial p_k} \otimes
          \frac{\partial}{\partial q^k}
      \right)}
    f \otimes g,
\end{equation}
where $\mu (f \otimes g) = fg$ is the pointwise product of 
$f, g \in C^\infty (T^*\mathbb R^n)[[\lambda]]$. For convenience we
consider compactly supported functions 
$f \in C^\infty_0 (T^*\mathbb R^n)[[\lambda]]$ and use the
Schr{\"o}dinger functional
\begin{equation}
    \label{SchroedingerFunct}
    \omega(f) = \int_{\mathbb R^n} \iota^* f \; d^nq,
\end{equation}
where $\iota: \mathbb R^n \hookrightarrow T^*\mathbb R^n$ is the zero
section and $d^nq$ the usual volume form. Clearly $\omega$ is a
classically positive linear functional and it turns out that $\omega$
is positive for $\starWeyl$, too, without the need of further
deformations. The Gel'fand ideals are given by
\begin{equation}
    \label{WeylGelfandIdeal}
    \mathcal J 
    = \{ f \in C^\infty_0 (T^*\mathbb R^n) \; | \; \iota^*f = 0 \} 
    \quad \textrm{and}\quad
    \boldsymbol{\mathcal J} 
    = \{ f \in C^\infty_0 (T^*\mathbb R^n)[[\lambda]] 
    \; | \; \iota^* Nf = 0 \},
\end{equation}
respectively. Here $N = \exp(\frac{\lambda}{2\im}\Delta)$ with 
$\Delta = \frac{\partial^2}{\partial q^k \partial p_k}$. Since $N$ is
a bijection it follows that in this example $\boldsymbol{\mathcal J}$
is a deformation of $\mathcal J$ as a subspace in the sense of
\cite[Def.~30]{BHW2000}. Moreover, the representation spaces are just
\begin{equation}
    \label{WaveFunctions}
    \mathcal H = C^\infty_0 (\mathbb R^n)
    \quad
    \textrm{and}
    \quad
    \boldsymbol{\mathcal H} 
    = \mathcal H[[\lambda]] 
    = C^\infty_0 (\mathbb R^n)[[\lambda]]
\end{equation}
with the usual $L^2$ inner product induced by $d^nq$. Thus in this
case we simply have 
$\boldsymbol{\mathcal H}_0 = \lambda \boldsymbol{\mathcal H}$.
Finally the $^*$-representations are now given by
\begin{equation}
    \label{SchroedingerRep}
    \varrho (f) \psi = \iota^*(f \pi^*\psi) = (\iota^*f) \psi
    \quad
    \textrm{and}
    \quad
    \boldsymbol{\varrho} (f) \psi = \iota^* N (f \starWeyl \pi^* \psi).
\end{equation}
Again, these considerations are easily transfered to the general case
of a cotangent bundle with the star product $\starWeyl$ as constructed
in \cite{BNW99a,BNW98a,BNPW98} and the results are literally the
same. We shall omit the rather obvious details.

We conclude that in the case of the Schr{\"o}dinger functional the
correct classical limit can be obtained by naively `setting 
$\lambda = 0$' everywhere. On the other hand, in the Wick case this
does not yield an acceptable result since then the classical
representation space would be 
$\mathbb C[[\cc y^1, \ldots, \cc y ^n]]$ with a highly degenerate
inner product as obtained from (\ref{FockInnerProd}) by `setting
$\lambda = 0$'. Thus one needs the more general construction of
dividing by $\boldsymbol{\mathcal H}_0$ in this case.

Physically speaking this can be interpreted as follows. The naive
classical limit by setting $\lambda = 0$ does not give meaningful
results for the `wave functions' themselves, i.e. for the elements in
the representation space $\boldsymbol{\mathcal H}$, but only for the
measurable quantities, i.e. the Hermitian inner products. Hence they
should be used to control the classical limit.

%
% acknowledgments
%

\section*{Acknowledgments}

First I would like to thank Henrique Bursztyn for various discussions
concerning this topic as well as for many helpful comments on the
manuscript. Moreover, I would like to thank Pierre Bieliavsky, Martin
Bordemann, Michel Cahen, Simone Gutt, and Cornelius Paufler for
valuable remarks.

%
% bibliography
%

\begin{small}

\end{small}

%\begin{small}
%\bibliographystyle{ewde}
%\bibliography{articles,books,preprints,misc}

\begin{thebibliography}{10}

\bibitem {BFLS84}
{\sc Basart, H., Flato, M., Lichnerowicz, A., Sternheimer, D.: }\newblock {\em
  Deformation Theory applied to Quantization and Statistical Mechanics}.
\newblock Lett. Math. Phys.  {\bf 8} (1984), 483--394.

\bibitem {BFFLS78}
{\sc Bayen, F., Flato, M., Fr{{\o}}nsdal, C., Lichnerowicz, A., Sternheimer,
  D.: }\newblock {\em Deformation Theory and Quantization}.
\newblock Ann. Phys.  {\bf 111} (1978), 61--151.

\bibitem {BBG98}
{\sc Bertelson, M., Bieliavsky, P., Gutt, S.: }\newblock {\em Parametrizing
  Equivalence Classes of Invariant Star Products}.
\newblock Lett. Math. Phys.  {\bf 46} (1998), 339--345.

\bibitem {BCG97}
{\sc Bertelson, M., Cahen, M., Gutt, S.: }\newblock {\em Equivalence of Star
  Products}.
\newblock Class. Quantum Grav.  {\bf 14} (1997), A93--A107.

\bibitem {BHW2000}
{\sc Bordemann, M., Herbig, H.-C., Waldmann, S.: }\newblock {\em BRST
  Cohomology and Phase Space Reduction in Deformation Quantization}.
\newblock Commun. Math. Phys.  {\bf 210} (2000), 107--144.

\bibitem {BNPW98}
{\sc Bordemann, M., Neumaier, N., Pflaum, M.~J., Waldmann, S.: }\newblock {\em
  On representations of star product algebras over cotangent spaces on
  Hermitian line bundles}.
\newblock Preprint Freiburg FR-THEP-98/24  {\bf math.QA/9811055} (November
  1998).

\bibitem {BNW98a}
{\sc Bordemann, M., Neumaier, N., Waldmann, S.: }\newblock {\em Homogeneous
  Fedosov Star Products on Cotangent Bundles I: Weyl and Standard Ordering with
  Differential Operator Representation}.
\newblock Commun. Math. Phys.  {\bf 198} (1998), 363--396.

\bibitem {BNW99a}
{\sc Bordemann, M., Neumaier, N., Waldmann, S.: }\newblock {\em Homogeneous
  Fedosov star products on cotangent bundles II: GNS representations, the WKB
  expansion, traces, and applications}.
\newblock J. Geom. Phys.  {\bf 29} (1999), 199--234.

\bibitem {BRW98a}
{\sc Bordemann, M., R{\"{o}}mer, H., Waldmann, S.: }\newblock {\em A Remark on
  Formal KMS States in Deformation Quantization}.
\newblock Lett. Math. Phys.  {\bf 45} (1998), 49--61.

\bibitem {BRW99}
{\sc Bordemann, M., R{\"o}mer, H., Waldmann, S.: }\newblock {\em KMS States and
  Star Product Quantization}.
\newblock Rep. Math. Phys.  {\bf 44} (1999), 45--52.

\bibitem {BW97a}
{\sc Bordemann, M., Waldmann, S.: }\newblock {\em A Fedosov Star Product of
  Wick Type for K{\"{a}}hler Manifolds}.
\newblock Lett. Math. Phys.  {\bf 41} (1997), 243--253.

\bibitem {BW97b}
{\sc Bordemann, M., Waldmann, S.: }\newblock {\em Formal GNS Construction and
  WKB Expansion in Deformation Quantization}.
\newblock In: {\sc Sternheimer, D., Rawnsley, J., Gutt, S. (eds.): }\newblock
  {\em Deformation Theory and Symplectic Geometry}, {\em Mathematical Physics
  Studies} no. {\bf 20},   315--319. Kluwer Academic Publisher, Dordrecht,
  Boston, London, 1997.

\bibitem {BW98a}
{\sc Bordemann, M., Waldmann, S.: }\newblock {\em Formal GNS Construction and
  States in Deformation Quantization}.
\newblock Commun. Math. Phys.  {\bf 195} (1998), 549--583.

\bibitem {BuWa99a}
{\sc Bursztyn, H., Waldmann, S.: }\newblock {\em Algebraic Rieffel Induction,
  Formal Morita Equivalence and Applications to Deformation Quantization}.
\newblock Preprint  {\bf math.QA/9912182} (December 1999).
\newblock To appear in J. Geom. Phys.

\bibitem {BuWa2000b}
{\sc Bursztyn, H., Waldmann, S.: }\newblock {\em {$^*$}-Ideals and Formal
  Morita Equivalence of {$^*$}-Algebras}.
\newblock Preprint  {\bf math.QA/0005227} (May 2000).
\newblock To appear in Int. J. Math.

\bibitem {BuWa2000c}
{\sc Bursztyn, H., Waldmann, S.: }\newblock {\em Deformation Quantization of
  Hermitian Vector Bundles}.
\newblock Preprint  {\bf math.QA/0009170} (September 2000).
\newblock To appear in Lett. Math. Phys.

\bibitem {BuWa2000a}
{\sc Bursztyn, H., Waldmann, S.: }\newblock {\em On Positive Deformations of
  {$^*$}-Algebras}.
\newblock In: {\sc Dito, G., Sternheimer, D. (eds.): }\newblock {\em
  Conf{\`e}rence Mosh{\`e} Flato 1999. Quantization, Deformations, and
  Symmetries}, {\em Mathematical Physics Studies} no. {\bf 22},   69--80.
  Kluwer Academic Publishers, Dordrecht, Boston, London, 2000.

\bibitem {CatFeld2000}
{\sc Cattaneo, A., Felder, G.: }\newblock {\em A Path Integral Approach to the
  Kontsevich Quantization Formula}.
\newblock Commun. Math. Phys.  {\bf 212} (2000), 591--611.

\bibitem {DL83b}
{\sc DeWilde, M., Lecomte, P. B.~A.: }\newblock {\em Existence of Star-Products
  and of Formal Deformations of the Poisson Lie Algebra of Arbitrary Symplectic
  Manifolds}.
\newblock Lett. Math. Phys.  {\bf 7} (1983), 487--496.

\bibitem {Dit90}
{\sc Dito, J.: }\newblock {\em Star-Product Approach to Quantum Field Theory:
  The Free Scalar Field}.
\newblock Lett. Math. Phys.  {\bf 20} (1990), 125--134.

\bibitem {DF2000}
{\sc D{\"u}tsch, M., Fredenhagen, K.: }\newblock {\em Algebraic Quantum Field
  Theory, Perturbation Theory, and the Loop Expansion}.
\newblock Preprint  {\bf hep-th/0001129} (January 2000).

\bibitem {Fed94a}
{\sc Fedosov, B.~V.: }\newblock {\em A Simple Geometrical Construction of
  Deformation Quantization}.
\newblock J. Diff. Geom.  {\bf 40} (1994), 213--238.

\bibitem {Fed96}
{\sc Fedosov, B.~V.: }\newblock {\em Deformation Quantization and Index
  Theory}.
\newblock Akademie Verlag, Berlin, 1996.

\bibitem {Fed98a}
{\sc Fedosov, B.~V.: }\newblock {\em Non-Abelian Reduction in Deformation
  Quantization}.
\newblock Lett. Math. Phys.  {\bf 43} (1998), 137--154.

\bibitem {Ger64}
{\sc Gerstenhaber, M.: }\newblock {\em On the Deformation of Rings and
  Algebras}.
\newblock Ann. Math.  {\bf 79} (1964), 59--103.

\bibitem {GS88}
{\sc Gerstenhaber, M., Schack, S.~D.: }\newblock {\em Algebraic Cohomology and
  Deformation Theory}.
\newblock In: {\sc Hazewinkel, M., Gerstenhaber, M. (eds.): }\newblock {\em
  Deformation Theory of Algebras and Structures and Applications},   13--264.
  Kluwer Academic Press, Dordrecht, 1988.

\bibitem {GR99}
{\sc Gutt, S., Rawnsley, J.: }\newblock {\em Equivalence of star products on a
  symplectic manifold; an introduction to Deligne's {\v{C}}ech cohomology
  classes}.
\newblock J. Geom. Phys.  {\bf 29} (1999), 347--392.

\bibitem {JS2000}
{\sc Jurco, B., Schupp, P.: }\newblock {\em Noncommutative Yang-Mills from
  equivalence of star products}.
\newblock Eur. Phys. J.  {\bf C14} (2000), 367--370.

\bibitem {JSW2000}
{\sc Jurco, B., Schupp, P., Wess, J.: }\newblock {\em Noncommutative gauge
  theory for Poisson manifolds}.
\newblock Nucl. Phys.  {\bf B584} (2000), 784--794.

\bibitem {Kon97b}
{\sc Kontsevich, M.: }\newblock {\em Deformation Quantization of Poisson
  Manifolds, I}.
\newblock Preprint  {\bf q-alg/9709040} (September 1997).

\bibitem {Kont99}
{\sc Kontsevich, M.: }\newblock {\em Operads and Motives in Deformation
  Quantization}.
\newblock Lett. Math. Phys.  {\bf 48} (1999), 35--72.

\bibitem {NT95a}
{\sc Nest, R., Tsygan, B.: }\newblock {\em Algebraic Index Theorem}.
\newblock Commun. Math. Phys.  {\bf 172} (1995), 223--262.

\bibitem {NT95b}
{\sc Nest, R., Tsygan, B.: }\newblock {\em Algebraic Index Theorem for
  Families}.
\newblock Adv. Math.  {\bf 113} (1995), 151--205.

\bibitem {OMY91}
{\sc Omori, H., Maeda, Y., Yoshioka, A.: }\newblock {\em Weyl Manifolds and
  Deformation Quantization}.
\newblock Adv. Math.  {\bf 85} (1991), 224--255.

\bibitem {Schomerus99}
{\sc Schomerus, V.: }\newblock {\em D-branes and deformation quantization}.
\newblock JHEP  {\bf 06} (1999).

\bibitem {Wal2000a}
{\sc Waldmann, S.: }\newblock {\em Locality in GNS Representations of
  Deformation Quantization}.
\newblock Commun. Math. Phys.  {\bf 210} (2000).

\end{thebibliography}
%\end{small}

\end{document}